\RequirePackage{ifpdf}
\ifpdf 
\documentclass[pdftex]{sigma}
\else
\documentclass{sigma}
\fi

\begin{document}
\allowdisplaybreaks

\renewcommand{\PaperNumber}{058}

\FirstPageHeading

\ShortArticleName{A Dual Mesh Method for a Non-Local Thermistor Problem}

\ArticleName{A Dual Mesh Method for a Non-Local\\ Thermistor Problem}

\Author{Abderrahmane EL HACHIMI~$^\dag$, Moulay Rchid SIDI
AMMI~$^\ddag$ and Delf\/im F.M. TORRES~$^\ddag$}
\AuthorNameForHeading{A. El Hachimi, M.R. Sidi Ammi and D.F.M.
Torres}

\Address{$^\dag$~UFR: Applied and Industrial Mathematics,
University of Chouaib Doukkali, El Jadida, Maroc}
\EmailD{\href{mailto:elhachimi@ucd.ac.ma}{elhachimi@ucd.ac.ma}}

\Address{$^\ddag$~Department of Mathematics, University of Aveiro,
3810-193 Aveiro, Portugal}
\EmailD{\href{mailto:sidiammi@mat.ua.pt}{sidiammi@mat.ua.pt},
\href{mailto:delfim@mat.ua.pt}{delfim@mat.ua.pt}}
\URLaddressD{\url{http://www.mat.ua.pt/delfim}}

\ArticleDates{Received December 20, 2005, in f\/inal form May 08,
2006; Published online June 02, 2006}

\Abstract{We use a dual mesh numerical method to study a non-local
parabolic problem arising from the well-known thermistor problem.}


\Classification{35K55; 65N15; 65N50}

\section{Introduction}

In this work we propose a dual mesh numerical scheme for analysis
of the following non-local parabolic problem coming from
conservation law of electric charges:
\begin{gather}
\frac{\partial u}{\partial t}-\nabla \cdot (k(u)\nabla u) =
\lambda \frac{f(u)}{ \left( \int_{\Omega} f(u)\, dx \right)^{2}}
\quad \mbox{in} \ \  \Omega \times ]0;T[  , \nonumber\\
u = 0 \quad \mbox{on} \ \  \partial \Omega \times ]0;T[  , \qquad
u/_{t=0}= u_0 \quad  \mbox{in} \ \  \Omega ,\label{11}
\end{gather}
where $\nabla$ denotes the gradient with respect to the
$x$-variables. The nonlinear problem \eqref{11} is obtained, under
some simplif\/icative conditions, by reducing the well-known
thermistor problem (cf., e.g., \cite{lac1,lac2,tza}), which
consists of the heat equation, with joule heating as a source, and
subject to current conservation:
\begin{gather}\label{14}
u_{t}= \nabla \cdot \left(k(u)\nabla u\right) + \sigma
(u)\left|\nabla \varphi\right|^{2}  , \qquad \nabla \cdot
\left(\sigma (u) \nabla \varphi\right) = 0 ,
\end{gather}
where the domain $\Omega \subset \mathbb{R}^{2}$ occupied by the
thermistor is a bounded convex polygonal; $\varphi = \varphi
(x,t)$ and $u = u(x,t)$ are, respectively, the distributions of
the electric potential and the temperature in $\Omega$;
$\sigma(u)$ and $k(u)$ are, respectively, the
temperature-dependant electrical and thermal conductivities;
$\sigma (u)\left|\nabla \varphi\right|^{2}$ is the joule heating.
The literature on problem \eqref{14} is vast (see
e.g.~\cite{ac,ci1,ci2,es1,es2,es3,es4,xu1,xu2}). With respect to
numerical approximation results to problem \eqref{14} we are aware
of \cite{aya,es4,cs,yx}: in~\cite{yx} a numerical analysis of the
non-steady thermistor problem by a f\/inite element method is
discussed; in \cite{cs} the authors study a spatially and
completely discrete f\/inite element model; in \cite{es4} a
semi-discretization by the backward Euler scheme is given for the
special case $k= Id$; in \cite{aya} a box approximation scheme is
presented and analyzed. A completely discrete scheme based on the
backward Euler method with semi-implicit linearization to
\eqref{14} is presented in \cite{cs} for the special case
$k(u)=1$. Existence and uniqueness of solutions to the problem
\eqref{11} were proved in \cite{es3}.

Finite volume methods emerged recently and seem to have a
signif\/icant role on concrete applications, because they have
very interesting properties in view of the subjacent physical
problems: in particular in conservation of f\/lows. An equation
coming from a conservation law has a good chance to be correctly
discretized by the f\/inite volume method. We also recall that
these schemes have been widely used to approximate solutions of
the heat linear equation, semi-linear or parabolic equations.
Since we consider data $f$ with lack of regularity when compared
to previous work, we need a new way to discretize \eqref{11}. We
present a dual mesh method capable of handling the non local term
$\frac{\lambda f(u)}{ \left( \int_{\Omega} f(u)\, dx \right)^{2}}$
which is a noticeable feature of \eqref{11}, by generalizing the
results of \cite{aya}. A box approximation scheme for discretizing
\eqref{11} with the case $k$ being dif\/ferent from the identity
is obtained. Speed of convergence is directly related with
regularity of the continuous problem. When one increases
regularity of the second term and data, the solution see its
regularity increasing in parallel, and precise speed of
convergence can be established. In the existing literature (see
e.g.~\cite{clt,cs}) the error estimates for both the f\/inite
element or volume element method are usually derived for solutions
that are suf\/f\/iciently smooth. Because the domain is polygonal,
special attention has to be paid to  regularity of the exact
solution. We give suf\/f\/icient conditions in terms of data and
the solution $u$ that yield error estimates (see hypothesis (H1)
below).

The text is organized as follows. In Section~\ref{sec2} we set up
the notation and the functional spaces used throughout the paper.
Section~\ref{sec3} introduces a box scheme model for problem
\eqref{11}, and  existence and uniqueness of the solution of the
approximating problem \eqref{31} is obtained from the f\/ixed
point theorem and  equivalence of norms in the f\/inite
dimensional space $S_{h}^{0}$. Finally, in Section~\ref{sec4},
under some regularity assumptions, we prove error estimates.

\section{Notation and functional spaces}
\label{sec2}

\begin{figure}[t]
\centerline{\includegraphics[width=5in]{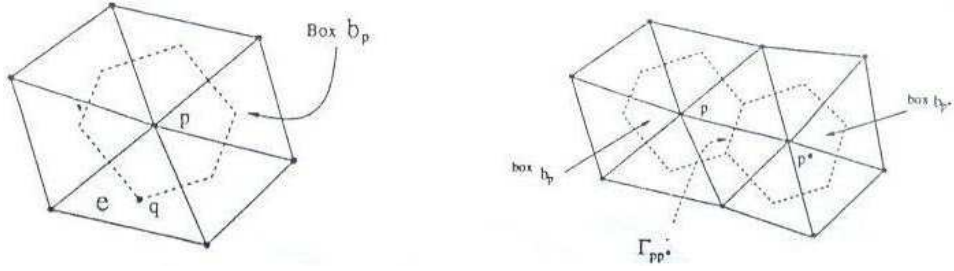}}
\caption{Construction of the dual mesh.} \label{fig1}
\end{figure}

Let $(\cdot,\cdot)$ and $\|\cdot\|$ denote the inner product and
norm in $L^{2}(\Omega)$; $H_{0}^{1}(\Omega)= \left\{ u \in
H^{1}(\Omega), u/\partial \Omega =0 \right\}$; $\|\cdot\|_{s}$,
$\|\cdot\|_{s,p}$ denote the $H^{s}(\Omega)$ and the
$W^{s,p}(\Omega)$ norm respectively; $T_{h}$ denote a
triangulation of~$\Omega$; $T_{v}^{h}$ be the set of vertices of a
quasi-uniform triangulation $T_{h}$; and
$\left\{S_{h}^{0}\right\}_{h>0}$ be the family of approximating
subspaces of $H_{0}^{1}(\Omega)$ def\/ined by
\begin{gather*}
S_{h}^{0} = \left\{v \in H_{0}^{1}(\Omega) : v/e \  \mbox{is a
linear function for all} \ e \in T_{h} \right\}  .
\end{gather*}
In the remainder of this paper we denote by $c$ various constants
that may depend on the data of the problem, and that are not
necessarily the same at each occurrence. We assume that the family
of triangulations is such that the following estimates \cite{cia}
hold for all $v \in S_{h}^{0}$:
\begin{gather}
\|v\|_{\beta , q} \leq ch^{r-\beta -2\max \{0,
1/p-1/q\}}\|v\|_{r,p} ,
\qquad 0\leq r \leq \beta \leq 1, \quad  1\leq p,q \leq \infty  , \nonumber\\
\|v\|_{0 , \infty} \leq c|\ln h|^{\frac{1}{2}} \|v\|_{1}
.\label{16}
\end{gather}
Let $P_{h}: L^{2}(\Omega)\rightarrow S_{h}^{0}$ be the standard
$L^{2}$-projection. One has \cite{cia}:
\begin{gather}
\|v- P_{h}v\|+h \|v- P_{h}v\|_{1} \leq ch^{2} \|v\|_{2} , \nonumber\\
\|v- P_{h}v\|_{0, \infty } \leq ch \|v\|_{2}, \qquad
\|P_{h}v\|_{1, \infty } \leq c \|v\|_{1, \infty } .\label{18}
\end{gather}
We construct the box scheme $B_{h}$ (dual mesh) employed in the
discretization as follows. From a~given triangle $e \in T_{h}$, we
choose a point $q \in \overline{e}$ as the intersection of the
perpendicular bisectors of the three edges of $e$. Then, we
connect $q$ by straight-line segments to the edge midpoints
of~$e$. To each vertex $p \in T_{v}^{h}$, we associate the box
$b_{p} \in B_{h}$, consisting of the union of subregions which
have $p$ as a corner (see Fig.~\ref{fig1}). For the piecewise
constant interpolation operator $I_{h}$, def\/ined by
\begin{gather*}
I_{h}: \mathcal{C}(\Omega)\rightarrow L^{2}(\Omega)  , \qquad
I_{h}v =v(p), \quad \mbox{on} \ \ b_{p} \in B_{h}, \quad \forall
\; p \in T_{v}^{h}  ,
\end{gather*}
we have the following standard error estimates \cite{aya,cai}:
\begin{gather}
c^{-1}\|v \| \leq \|I_{h}v\| \leq c \|v\| ,
\quad \forall \; v \in S_{0}^{h} , \nonumber\\
\| v-I_{h}v\| \leq c h\|v\|_{1}, \quad \forall \; v \in S_{0}^{h}
.\label{110}
\end{gather}
We denote by $N_{h}(p)$ the set of the neighboring vertices of $p
\in T_{v}^{h}$, $\partial b = \bigsqcup_{p \in T_{v}^{h}}\partial
b_{p}$, $\partial b_{p}= \bigsqcup_{p*\in N_{h}(p)}
\{\Gamma_{pp*}\}$, where $\Gamma_{pp*}= \partial b_{p} \bigcap
\partial b_{p*}$ (see Fig.~\ref{fig1}). Let $l_{\partial b}:
\partial b \rightarrow \mathbb{R}^{+}$ be def\/ined as follows:
for $p\in T_{v}^{h}$ and $b_{p} \in B_{h}$,
\begin{gather*}
l_{\partial b}/\Gamma_{pp*}= |p-p*| \quad \mbox{for} \ \ p*\in
N_{h}(p) .
\end{gather*}
For $b \in B_{h}$, we denote the jump in $w$ across $\partial b$
at $x$ by $[w]_{\partial b}(x)= w(x+0)-w(x-0)$, where $w(x\pm 0)$
are the outside and inside limit values of $w(x)$ along the normal
directions for $\partial b$.

We now collect from the literature \cite{aya,cai} some important
lemmas and trace results, that are needed in the sequel.
\begin{lemma}
\label{lm21} Assume that $B_{h}$ is a dual mesh. If $v$ is a
piecewise linear function, and $x$ is not a~vertex, then
\begin{gather*}
[I_{h}v]/_{\partial b_{p}}(x)= \frac{\partial v}{\partial n}
l_{\partial b}/\Gamma_{pp*}, \quad x \in \Gamma_{pp*}, \quad
\forall \; b\in B_{h} ,
\end{gather*}
where $n$ is the unit outward normal vector on $\partial b$.
\end{lemma}

The $h$-dependent norms are def\/ined as follows:
\begin{gather*}
\|v\|_{1, h} =\left(\sum_{l \in \partial
b}\left|[I_{h}v]_{l}\right|^{2}\right)^{\frac{1}{2}} \quad
\mbox{and} \quad  \|v\|_{0, h} = \left\|I_{h}v\right\| .
\end{gather*}

\begin{lemma}
\label{lm22} There exists a constant $c > 0$ such that
\begin{gather*}
c^{-1}\left\|\nabla v\right\| \leq \|v\|_{1,h} \leq c \|\nabla v\|
,
\quad \forall\; v \in S_{h}^{0}  , \\
c^{-1} \| v\| \leq \|v\|_{0,h} \leq c \| v \|  , \quad \forall \;
v \in S_{h}^{0}  .
\end{gather*}
\end{lemma}

\begin{lemma}
\label{lm23} For any $a \in \mathcal{C}(\overline{\Omega})$ there
exists a positive constant $c$ such that
\begin{gather}
\label{112} \left|- \sum_{b \in B_{h}}\int_{\partial b}a
\frac{\partial u}{\partial n} I_{h}v\right| \leq
c\|u\|_{1}\|v\|_{1}, \quad \forall \; u, v \in S_{h}^{0} .
\end{gather}
Moreover, if there exists a constant $a_{0}>0$ such that $a\geq
a_{0}$ in $\Omega$, then
\begin{gather}
\label{113} c^{-1}\|v\|_{1}^{2} \leq - \sum_{b \in
B_{h}}\int_{\partial b}a \frac{\partial v}{\partial n} I_{h}v ,
\quad \forall\;  v \in S_{h}^{0} .
\end{gather}
\end{lemma}

Let $Q_{h}: H^{2}(\Omega)\rightarrow S_{h}^{0}$ be def\/ined by
$Q_{h}u-i_{h}u \in S_{h}^{0}$, and
\begin{gather}
\label{114} - \sum_{b \in B_{h}}\int_{\partial b}a \frac{\partial
(u-Q_{h}u)}{\partial n} I_{h}v =0, \quad \forall \; v \in
S_{h}^{0} ,
\end{gather}
where $i_{h}: \mathcal{C}(\Omega)\rightarrow S_{h}^{0}$ is the
Lagrangian interpolation operator and $u \in H^{2}(\Omega)$.

\begin{lemma}
\label{lm24} Assume that $a \in L^{\infty}(\Omega)$, with $a \geq
a_{0}$ for some constant $a_{0}>0$. Then, there exists $c>0$ such
that for $u\in  H^{2}(\Omega)$
\begin{gather}
\label{115} \|u-Q_{h}u\|_{1} \leq ch \|u\|_{2} .
\end{gather}
Moreover, if $u \in H^{2}(\Omega)\bigcap W^{1, \infty}(\Omega)$,
then
\begin{gather}
\label{116} \|Q_{h}u\|_{1,\infty} \leq c\left(\|u\|_{1,
\infty}+\|u\|_{2}\right) .
\end{gather}
\end{lemma}

\begin{lemma}
\label{lm25} For each $b\in B_{h}$ one has
\begin{gather*}
h^{\frac{1}{2}}\left\|v\right\|_{L^{2}\left(\partial b\right)}
\leq c \left( \|v\|_{L^{2}(b)} + h \|v\|_{H^{1}(b)}\right), \quad
\forall \; v \in H^{1}(b)  .
\end{gather*}
\end{lemma}

Throughout this work, we assume that the following hypotheses on
the solution and data of problem \eqref{11} are satisf\/ied:
\begin{itemize}\itemsep=0pt

\item[(H1)] $u \in L^{\infty}(H_{0}^{1}(\Omega)\bigcap
H^{2}(\Omega))$, $u_{t}\in L^{2}(H^{1}(\Omega))$;

\item[(H2)] $c^{-1}\leq k(s)\leq c $;

\item[(H3)] there exist positive constants  $c_1$, $c_2$ and
$\nu$, such that $\nu \leq f(\xi )\leq c_{1}| \xi|+c_{2}$ for all
$\xi \in \mathbb{R}$;

\item[(H4)] $\left|f(\xi)-f(\xi')\right| +
\left|k(\xi)-k(\xi')\right| \leq c \left|\xi - \xi'\right|$.

\end{itemize}

\section{Existence and uniqueness result for the box scheme method}
\label{sec3}

Let $u$ be the solution of \eqref{11}. Integrating over an element
$b$ in $B_{h}$ we obtain:
\begin{gather}
\label{31a} \int_{b}u_{t} - \int_{\partial b} k(u)\frac{\partial
u}{\partial n} = \frac{\lambda }{ \left(\int f(u)\, dx
\right)^{2}} \int_{b}f(u)  , \quad \forall \; b \in B_{h} .
\end{gather}
We consider a box scheme def\/ined as follows: f\/ind $u_{h}\in
S_{h}^{0}$ such that
\begin{gather}
\label{31} \big(I_{h}u_{t}^{h}, I_{h}v\big) - \sum_{b \in
B_{h}}\int_{\partial b}k\big(u^{h}\big) \frac{\partial
u^{h}}{\partial n} I_{h}v = \frac{\lambda}{\left(\int_{\Omega}
f(u^{h})\, dx \right)^{2}} \big(f\big(u^{h}\big), I_{h}v\big) ,
\quad \forall \;  v \in S_{h}^{0}  ,
\end{gather}
where $u^{h}(0)=P_{h}u_{0}$ and $I_{h}$ is the interpolation
operator.

\begin{theorem}
\label{thm31} Let {\rm (H1)--(H4)} be satisfied. Then, for each
$h>0$, there exists $t_{0}(h)$ such that \eqref{31} possesses  a
unique solution $u^{h}$ for $0\leq t \leq t_{0}(h)$.
\end{theorem}
\begin{proof}
We begin by proving existence of solution. We def\/ine a nonlinear
operator $G$ from~$S_{h}^{0}$ to~$S_{h}^{0}$ as follows. For each
$u^{h} \in S_{h}^{0}$, $w^{h}= G(u^{h})$ is obtained as the unique
solution of the following problem:
\begin{gather} \label{32}
\big(I_{h}w_{t}^{h}, I_{h}v\big)- \sum_{b \in B_{h}}\int_{\partial
b}k\big(u^{h}\big) \frac{\partial w^{h}}{\partial n} I_{h}v =
\frac{\lambda}{\left( \int_{\Omega} f(u^{h})\, dx \right)^{2}}
\big(f\big(u^{h}\big), I_{h}v\big), \quad \forall \; v \in
 S_{h}^{0}  .
\end{gather}
We remark that $G$ is well def\/ined. Using $v=w^{h}$ as a test
function in \eqref{32}, hypotheses (H2) and~(H3), and Holder's
inequality, we can write:
\begin{gather*}
\frac{1}{2}\frac{d}{dt}\|I_{h}w^{h}\|^{2}+c \|w^{h}\|_{1}^{2} \leq
c \big(f\big(u^{h}\big), I_{h}w^{h}\big)
 \leq c\int \big(|u^{h}|+1\big) |I_{h}w^{h}| \\
\qquad{}\leq c \|u^{h}\|_{L^{2}}\|I_{h}w^{h}\|_{L^{2}} + c
\|I_{h}w^{h}\|
 \leq c \|u^{h}\|_{1}\|I_{h}w^{h}\|_{1}+ c \|I_{h}w^{h}\|\\
\qquad{}
 \leq \frac{c}{2} \|I_{h}w^{h}\|_{1}^{2}+c \|u^{h}\|_{1}^{2} + c.
\end{gather*}
Thus, we have
\begin{gather}
\label{33} \frac{d}{dt}\|I_{h}w^{h}\|^{2}+c \|w^{h}\|_{1}^{2}
  \leq c \|u^{h}\|_{1}^{2} + c \, .
\end{gather}
Integrating \eqref{33} with respect to $t$ and using the
equivalency of $\|I_{h}\cdot\|$ and $\|\cdot\|$ in $S_{h}^{0}$
(see \eqref{110}) yields
\begin{gather*}
 \|w^{h}\|^{2}+c \int_{0}^{t}\|w^{h}\|^{2}_{1} \leq c
\|I_{h}P_{h}u_{0}\|^{2} + c \int_{0}^{t} \|u^{h}\|_{1}^{2} \, dx +
ct\\
\qquad {} \leq c \|u_{0}\|^{2} + c \int_{0}^{t} \|u^{h}\|_{1}^{2}
\, dx + ct .
\end{gather*}
Def\/ine now the following set
\begin{gather*}
D = \left\{ u^{h}\in S_{h}^{0}, \|u^{h}\|^{2}+c
\int_{0}^{t}\|u^{h}\|^{2}_{1}\leq c \left(\|u_{0}\|^{2}+ 1\right)
\right\}  .
\end{gather*}
We can easily see that $D$ is closed subset of $L^{\infty}(0, t,
L^{2}(\Omega))$ with its natural norm.
 We conclude that there exists $t>0$
such that $G(D) \subset D$.
 To obtain that $G$ has a f\/ixed point $w^{h}=G(w^{h})$,
we prove that $G$ is a contraction. Conclusion follows from
Banach's f\/ixed point theorem. For this purpose, let $u_{1}^{h}$
and $u_{2}^{h}\in S_{h}^{0}\times S_{h}^{0}$ such that
$Gu_{1}^{h}= w_{1}^{h}$ and $Gu_{2}^{h}= w_{2}^{h}$. We have, from
the equation \eqref{32} verif\/ied by $w_{1}^{h}$ and $w_{2}^{h}$,
that
\begin{gather*}
\big(I_{h}\big(w_{1t}^{h}-w_{2t}^{h}\big), I_{h}v\big) - \sum_{b
\in B_{h}}\int_{\partial b}k\big(u^{h}_{1}\big) \frac{\partial
w^{h}_{1}}{\partial n} I_{h}v + \sum_{b \in B_{h}}\int_{\partial
b}k\big(u^{h}_{2}\big) \frac{\partial w^{h}_{2}}{\partial n} I_{h}v \\
\qquad{} = \frac{\lambda}{\left( \int_{\Omega} f(u^{h}_{1})\, dx
\right)^{2}} \big(f\big(u^{h}_{1}\big), I_{h}v\big) -
\frac{\lambda}{\left( \int_{\Omega} f(u^{h}_{2})\, dx \right)^{2}}
\big(f\big(u^{h}_{2}\big), I_{h}v\big)  .
\end{gather*}
On the other hand, one has
\begin{gather*}
 -\sum_{b \in B_{h}}\int_{\partial b}k\big(u^{h}_{1}\big) \frac{\partial
w^{h}_{1}}{\partial n} I_{h}v + \sum_{b \in B_{h}}\int_{\partial
b}k\big(u^{h}_{2}\big) \frac{\partial w^{h}_{2}}{\partial n} I_{h}v \\
\qquad{} = -\sum_{b \in B_{h}}\int_{\partial b}k(u^{h}_{1})
\frac{\partial \big(w^{h}_{1}-w^{h}_{2}\big)}{\partial n} I_{h}v
+\sum_{b \in B_{h}}\int_{\partial
b}\big(k\big(u^{h}_{2}\big)-k\big(u^{h}_{1}\big)\big)
\frac{\partial w^{h}_{2}}{\partial n} I_{h}v  .
\end{gather*}
Schwartz inequality implies that
\begin{gather*}
\sum_{b \in B_{h}}\int_{\partial
b}\big(k\big(u^{h}_{1}\big)-k\big(u^{h}_{2}\big)\big)
\frac{\partial w^{h}_{2}}{\partial n} I_{h}v \geq -c
\|w_{2}^{h}\|_{1, \infty }\left(\sum_{b \in
B_{h}}\|u_{1}^{h}-u_{2}^{h}\|_{0,\partial b}\right)\|v\|  .
\end{gather*}
By Lemma \ref{lm25}, we have
\begin{gather*}
  h^{\frac{1}{2}}\sum_{b
\in B_{h}}\|u_{1}^{h}-u_{2}^{h}\|_{0,\partial b} \leq c \sum_{b
\in B_{h}}\big(\|u_{1}^{h}-u_{2}^{h}\|_{L^{2}(b)}+h
\|u_{1}^{h}-u_{2}^{h}\|_{H^{1}(b)}\big) \leq c
\|u_{1}^{h}-u_{2}^{h}\|_{1} .
\end{gather*}
Thus, from the inverse estimate \eqref{16},
\begin{gather} \label{35}
\sum_{b \in B_{h}}\int_{\partial
b}\big(k\big(u^{h}_{1}\big)-k\big(u^{h}_{2}\big)\big)
\frac{\partial w^{h}_{2}}{\partial n} I_{h}v \geq
-c(h)\|w_{2}^{h}\|_{1}\|u_{1}^{h}-u_{2}^{h}\|_{1} \|v\|_{1} .
\end{gather}
On the basis of hypotheses (H1)--(H4), we have:
\begin{gather}
\frac{\lambda}{\left( \int_{\Omega} f(u^{h}_{1})\, dx \right)^{2}}
\big(f\big(u^{h}_{1}\big), I_{h}v\big)- \frac{\lambda}{\left(
\int_{\Omega}
f\left(u^{h}_{2}\right)\, dx \right)^{2}} \big(f\big(u^{h}_{2}\big), I_{h}v\big)\nonumber\\
\qquad{}= \frac{\lambda}{\left( \int_{\Omega}
f\left(u^{h}_{1}\right)\, dx \right)^{2}}
\big(f\big(u^{h}_{1}\big)-f\big(u^{h}_{2}\big), I_{h}v\big) \nonumber\\
\qquad{} + \lambda \left(\frac{1}{\left( \int_{\Omega}
f\left(u^{h}_{1}\right)\, dx \right)^{2}}-\frac{1}{\left(
\int_{\Omega} f\left(u^{h}_{2}\right)\, dx
\right)^{2}}\right)\big(f\big(u^{h}_{2}\big),
I_{h}v\big)\nonumber\\
\qquad {}\leq c \|u_{1}^{h}-u_{2}^{h}\| \, \|v \| +
\lambda\frac{\left(\int_{\Omega}
f\left(u^{h}_{2}\right)-f\left(u^{h}_{1}\right)\right)\left(\int_{\Omega}
f\left(u^{h}_{2}\right)+f\left(u^{h}_{1}\right)\right)}{\left(\int_{\Omega}
f\left(u^{h}_{2}\right)\, dx \right)^{2}\left(\int_{\Omega}
f\left(u^{h}_{1}\right)\, dx
\right)^{2}}\big(f\big(u^{h}_{2}\big),
I_{h}v\big)\nonumber\\
\qquad{} \leq c \|u_{1}^{h}-u_{2}^{h}\| \, \|v \| + c
\|I_{h}v\|_{L^{2}(\Omega)}\|u_{1}^{h}-u_{2}^{h}\|_{L^{1}(\Omega)}\nonumber\\
\qquad{} \leq c \|u_{1}^{h}-u_{2}^{h}\| \, \|v \| \leq c
\|u_{1}^{h}-u_{2}^{h}\|_{1} \, \|v \|_{1}  .\label{36}
\end{gather}
It follows from \eqref{35} and \eqref{36} that
\begin{gather}
\big(I_{h}\big(w_{1t}^{h}-w_{2t}^{h}\big), I_{h}v\big) - \sum_{b
\in B_{h}}\int_{\partial b}k(u^{h}_{1}) \frac{\partial
\big(w^{h}_{1}-w^{h}_{2}\big)}{\partial n} I_{h}v  \nonumber\\
\leq c(h)\|u_{1}^{h}-u_{2}^{h}\|_{1} \|v\|_{1} -\sum_{b \in
B_{h}}\int_{\partial
b}\big(k\big(u^{h}_{1}\big)-k\big(u^{h}_{2}\big)\big)
\frac{\partial w^{h}_{2}}{\partial n} I_{h}v \leq c(h)
\|u_{1}^{h}-u_{2}^{h}\|_{1} \|v \|_{1}  .\label{38}
\end{gather}
Now, using $v=w_{1}^{h}-w_{2}^{h}$ as a test function in
\eqref{38}, we obtain from \eqref{113}:
\begin{gather}
\label{377}
\frac{1}{2}\frac{d}{dt}\big\|I_{h}\big(w_{1}^{h}-w_{2}^{h}\big)\big\|^{2}
+ c \|w_{1}^{h}-w_{2}^{h}\|^{2}_{1} \leq c(h)
\|u_{1}^{h}-u_{2}^{h}\|_{1}\|w_{1}^{h}-w_{2}^{h}\|_{1}  .
\end{gather}
With use of the Holder's inequality and  equivalency of
$\|I_{h}\cdot\|$ and $\|\cdot\|$,  integration of \eqref{377} with
respect to time gives:
\begin{gather*}
\big\|\big(w_{1}^{h}-w_{2}^{h}\big)\big\|^{2} \leq c
\big\|I_{h}\big(w_{1}^{h}-w_{2}^{h}\big)\big\|^{2} \leq c(h)
\int_{0}^{t} \big\|\big(u_{1}^{h}-u_{2}^{h}\big)\big\|^{2}_{1}ds .
\end{gather*}
Thus $G$ is a contraction. We prove now uniqueness. Following the
same arguments as before, we have
\begin{gather}
\label{39} \big(I_{h}\big(u_{1t}^{h}-u_{2t}^{h}\big), I_{h}v\big)
- \sum_{b \in B_{h}}\int_{\partial b}k(u^{h}_{1}) \frac{\partial
(u^{h}_{1}-u^{h}_{2})}{\partial n} I_{h}v \leq c(h)
\|u_{1}^{h}-u_{2}^{h}\|_{1} \|v \|  .
\end{gather}
Choosing $v= u_{1}^{h}-u_{2}^{h}$ as test function in \eqref{39},
using again \eqref{113} and integrating, we obtain
\begin{gather*}
\big\|\big(u_{1}^{h}-u_{2}^{h}\big)\big\|^{2} \leq
c(h)\int_{0}^{t}\big\|\big(u_{1}^{h}-u_{2}^{h}\big)\big\|^{2}ds  ,
\end{gather*}
which gives, by Gronwall's Lemma,  uniqueness of solution.
\end{proof}

\section{Error analysis}
\label{sec4}

In this section we prove error estimates under certain assumptions
on  regularity of the exact solution $u$.

\begin{theorem}
\label{thm2} Under assumptions {\rm (H1)--(H4)}, if
$\left(u,u^{h}\right)$ are solutions of \eqref{31a}--\eqref{31}
for $0\leq t \leq t_{0}(h)$, then
\[
\|u^{h}-u\|_{L^{\infty}(L^{2})} + \|u^{h}-u\|_{L^{2}(H^{1})} \leq
ch .
\]
\end{theorem}

\begin{proof}
From \eqref{31a} and \eqref{31} we obtain
\begin{gather}
\big(I_{h}\big(u^{h}-P_{h}u\big)_{t}, I_{h}v\big) - \sum_{b \in
B_{h}}\int_{\partial b}k\big(u^{h}\big)
\frac{\partial \left(u^{h}_{1}-P_{h}u\right)}{\partial n} I_{h}v \nonumber\\
= \frac{\lambda}{\left(\int_{\Omega} f\left(u^{h}\right)\, dx
\right)^{2}} \big(f\big(u^{h}\big), I_{h}v\big) -
\frac{\lambda}{\left(\int_{\Omega} f(u)\, dx\right)^{2}}
\big(f(u), I_{h}v\big) + \!\sum_{b \in B_{h}}\!\int_{\partial b}
\! k\big(u^{h}\big)
\frac{\partial \left(P_{h}u-u\right)}{\partial n} I_{h}v \nonumber\\
-\sum_{b \in B_{h}}\int_{\partial
b}\big(k(u)-k\big(u^{h}\big)\big) \frac{\partial u}{\partial n}
I_{h}v + \big((I-P_{h})u_{t},I_{h}v\big) +
\big((I-I_{h})P_{h}u_{t}, I_{h}v\big)  .\label{41}
\end{gather}
We now estimate, separately, the terms on the right-hand side of
\eqref{41}. We have from \eqref{112} and~\eqref{18} that
\begin{gather} \label{42}
\left|\sum_{b \in B_{h}}\int_{\partial b}k\big(u^{h}\big)
\frac{\partial (P_{h}u-u)}{\partial n} I_{h}v\right| \leq c
\|P_{h}u-u\|_{1}\|v\|_{1}\leq ch\|u\|_{2}\|v\|_{1}  \leq ch
\|v\|_{1} ,
\\
\left|\sum_{b \in B_{h}}\int_{\partial
b}\big(k(u)-k\big(u^{h}\big)\big)
\frac{\partial u}{\partial n} I_{h}v\right|\nonumber\\
\qquad{} \leq c \|v\|_{1}\left(\sum_{b\in
B_{h}}\left(\int_{\partial b}|u-u^{h}|\left|\frac{\partial
u}{\partial n}\right|\right)^{2}\right)^{\frac{1}{2}} \leq c
h^{\frac{1}{2}}\|v\|_{1}\left(\sum_{b \in
B_{h}}\|u-u^{h}\|^{2}_{0,\partial b}\right)^{\frac{1}{2}}  .
\label{43}
\end{gather}
By Lemma~\ref{lm25}, inverse inequality \eqref{16} and \eqref{18},
we have:
\begin{gather}
\sum_{b \in B_{h}} \|u^{h}-u\|_{0, \partial b}^{2} \leq 2\sum_{b
\in B_{h}} \|u^{h}-P_{h}u\|_{0, \partial b}^{2}+2\sum_{b \in
B_{h}} \|P_{h}u-u\|_{0, \partial b}^{2}\nonumber\\
\phantom{\sum_{b \in B_{h}} \|u^{h}-u\|_{0, \partial b}^{2}}{}
\leq c \big(h^{2}+ \|u^{h}-P_{h}u\|^{2}\big)  . \label{44}
\end{gather}
Consequently, we obtain from \eqref{44} that
\begin{gather} \label{45}
 \left| \sum_{b \in B_{h}}\int_{\partial b}\big(k(u)-k\big(u^{h}\big)\big)
\frac{\partial u}{\partial n} I_{h}v\right| \leq c
\big(h+\|u^{h}-P_{h}u\|\big)\|v\|_{1} .
\end{gather}
Based on our earlier development in \eqref{36}, we also know:
\begin{gather} \label{46}
\frac{\lambda}{\left( \int_{\Omega} f(u^{h})\, dx \right)^{2}}
\big(f\big(u^{h}\big), I_{h}v\big)- \frac{\lambda}{\left(
\int_{\Omega} f(u)\, dx \right)^{2}} (f(u), I_{h}v) \leq c
\|u^{h}-u\|\|v\|_{1} .
\end{gather}
Let $v=u^{h}-P_{h}u$ be a test function in \eqref{41}. Using
Lemma~\ref{lm23}, it follows from \eqref{42}--\eqref{46} that
\begin{gather*}
\frac{1}{2}\frac{d}{dt}  \|I_{h}(u^{h}-P_{h}u)\|^{2}+c
\|u^{h}-P_{h}u\|^{2}_{1} \\
\qquad{} \leq c\|u^{h}-u\|\|u^{h}-P_{h}u\|_{1} +
ch\|u^{h}-P_{h}u\|_{1}
+c\big(h+\|u^{h}-P_{h}u\|\big)\|u^{h}-P_{h}u\|_{1}\\
 \qquad\quad{} + c(\|(I-P_{h})u_{t}\|
+ \|(I-I_{h})P_{h}u_{t}\|)\|u^{h}-P_{h}u\|\\
 \qquad{}\leq c
(\|(I-P_{h})u_{t}\|
+ \|(I-I_{h})P_{h}u_{t}\|)\|u^{h}-P_{h}u\|\\
 \qquad \quad{} +c\big(h+\|u^{h}-P_{h}u\|\big)\|u^{h}-P_{h}u\|_{1}
+ \|P_{h}u-u\|\|u^{h}-P_{h}u\|_{1} .
\end{gather*}
By properties \eqref{18} and Cauchy's inequality, it follows:
\begin{gather*}
 \frac{d}{dt}\big\|I_{h}\big(u^{h}-P_{h}u\big)\big\|^{2}+c
\|u^{h}-P_{h}u\|^{2}_{1} \leq c\|P_{h} u-u\|_{1}\|u^{h}-P_{h}u\|_{1}\\
 \qquad \quad{} + c\big\{
h+\|(I-P_{h})u_{t}\|+\|(I-I_{h})P_{h}u_{t}\|+\|u^{h}-P_{h}u\|
\big\} \|u^{h}-P_{h}u\|_{1} \\
\qquad{} \le c \big\{
h^{2}+\|(I-P_{h})u_{t}\|^{2}_{1}+\|(I-I_{h})P_{h}u_{t}\|^{2}+\|u^{h}-P_{h}u\|^{2}
\big\}+\frac{c}{2} \|u^{h}-P_{h}u\|_{1}^{2}\\
\qquad{}\le c \big\{
h^{2}+h^{2}\|u_{t}\|_{2}^{2}+ch^{2}\|P_{h}u_{t}\|_{1}^{2}
\big\}+c\|u^{h}-P_{h}u\|^{2}+\frac{c}{2}\|u^{h}-P_{h}u\|_{1}^{2}.
\end{gather*}
Hence,
\begin{gather}\label{47}
 \frac{d}{dt}\big\|I_{h}\big(u^{h}-P_{h}u\big)\big\|^{2} + c
\big\|u^{h}-P_{h}u\big\|_{1}^{2} \leq ch^{2}+c\|u^{h}-P_{h}u\|^{2}
.
\end{gather}
Integrating \eqref{47} and applying Gronwall Lemma and using again
the equivalency of $\|\cdot\|$ and $\|I_{h}\cdot\|$, we get that
\[
\|u^{h}-P_{h}u\|^{2}+c \int_{0}^{t}\|u^{h}-P_{h}u\|_{1}^{2} \leq
ch^{2}  .
\]
Then, by the triangular inequality, we conclude with the intended
result.
\end{proof}

Under more restrictive hypotheses on the data, it is possible to
derive the following error estimate.

\begin{theorem}
Assume {\rm (H1)--(H4)}. If $k(s)=1$ and $u_{0}\in
H_{0}^{1}(\Omega)\bigcap H^{2}(\Omega )$, then
\begin{gather*}
\|u^{h}-u\|_{L^{\infty}(H^{1})}\leq ch .
\end{gather*}
\end{theorem}

\begin{proof}
From equations \eqref{31a} and \eqref{31}, we have:
\begin{gather*}
\big(I_{h}u_{t}^{h}-u_{t}, I_{h}v\big) - \sum_{b \in
B_{h}}\int_{\partial b} \frac{\partial u^{h}}{\partial n} I_{h}v+
\sum_{b \in B_{h}}\int_{\partial b}
\frac{\partial u}{\partial n}I_{h}v  \\
\qquad{}= \frac{\lambda}{\left( \int_{\Omega} f(u^{h})\, dx
\right)^{2}} \big(f\big(u^{h}\big), I_{h}v\big)-
\frac{\lambda}{\left( \int_{\Omega} f(u)\, dx \right)^{2}} (f(u),
I_{h}v) .
\end{gather*}
Using the def\/inition \eqref{114} of $Q_{h}$, we get
\begin{gather*}
\big(I_{h}u_{t}^{h}-u_{t}, I_{h}v\big) - \sum_{b \in
B_{h}}\int_{\partial b}
\frac{\partial \left(u^{h}-Q_{h}u\right)}{\partial n} I_{h}v \\
\qquad{}= \frac{\lambda}{\left( \int_{\Omega} f(u^{h})\, dx
\right)^{2}} \big(f\big(u^{h}\big), I_{h}v\big)-
\frac{\lambda}{\left( \int_{\Omega} f(u)\, dx \right)^{2}} (f(u),
I_{h}v) ,
\end{gather*}
and it follows that
\begin{gather}
\big(I_{h}\big(u_{t}^{h}-Q_{h}u\big)_{t}, I_{h}v\big) - \sum_{b
\in B_{h}}\int_{\partial b}
\frac{\partial (u^{h}-Q_{h}u)}{\partial n} I_{h}v \nonumber\\
\qquad{}= \frac{\lambda}{\left( \int_{\Omega} f(u^{h})\, dx
\right)^{2}} \big(f\big(u^{h}\big), I_{h}v\big) -
\frac{\lambda}{\left(
\int_{\Omega} f(u)\, dx \right)^{2}} (f(u), I_{h}v) \nonumber\\
\qquad\quad{}+((I-Q_{h})u_{t}, I_{h}v) + ((I-I_{h})Q_{h}u_{t},
I_{h}v) .\label{410}
\end{gather}
In order to estimate the right hand side of the last inequality,
we treat both terms separately. By similar arguments to those used
in \eqref{36},
\begin{gather*}
\left|\frac{\lambda}{\left( \int_{\Omega} f(u^{h})\, dx
\right)^{2}} \big(f\big(u^{h}\big), I_{h}v\big) -
\frac{\lambda}{\left(\int_{\Omega} f(u)\, dx \right)^{2}}
\left(f(u), I_{h}v\right)\right|\leq c\|u^{h}-u\| \|v\| .
\end{gather*}
Taking a function test $v=(u^{h}-Q_{h}u)_{t}$ in \eqref{410}, by
\eqref{115} and \eqref{116} we have
\begin{gather}
\|I_{h}(u^{h}  -Q_{h}u)_{t}\|^{2} - \sum_{b \in
B_{h}}\int_{\partial b} \frac{\partial (u^{h}-Q_{h}u)}{\partial n}
I_{h}\big(u^{h}-Q_{h}u\big)_{t}\nonumber \\
\qquad{} \leq c \big\{ h+ \|(I-I_{h})Q_{h}u_{t}\|+
\|(I-Q_{h})u_{t}\|+
\|Q_{h}u-u^{h }\|\big\} \|I_{h}(u^{h}-Q_{h}u)_{t}\|\nonumber\\
\qquad{} \leq c\big\{ h + ch \|u_{t}\|_{2}+ch\|Q_{h}u_{t}\|_{1}+
\|Q_{h}u-u^{h }\| \big\} \|I_{h}(u^{h}-Q_{h}u)_{t}\|\nonumber\\
\qquad{} \leq c\big\{ h +  \|Q_{h}u-u^{h }\| \big\} \|I_{h}(u^{h}-Q_{h}u)_{t}\|\nonumber\\
\qquad{} \leq
ch^{2}+c\|Q_{h}u-u^{h}\|^{2}+\frac{1}{2}\|I_{h}(u^{h}-Q_{h}u)_{t}\|^{2}
.\label{411}
\end{gather}
 Integrating \eqref{411}, we arrive to
\begin{gather*}
\|u^{h}-Q_{h}u\|^{2}_{1} \leq
c\left(h^{2}+\int_{0}^{t}\|u^{h}-Q_{h}u\|^{2}\right)   \\
\phantom{\|u^{h}-Q_{h}u\|^{2}_{1}}{} \leq
c\left(h^{2}+\int_{0}^{t}\|u^{h}-Q_{h}u\|_{1}^{2}\right) =  c
h^{2} + c \|u^{h}-Q_{h}u\|^{2}_{L^{2}(H^{1}(\Omega))}  ,
\end{gather*}
and Theorem~\ref{thm2} gives
\[
\|u^{h}-Q_{h}u\|^{2}_{1} \leq c h^{2}.
\]
On the other hand, by triangular inequality, \eqref{115} and the
regularity of the exact solution $u$, we have
\[
\|u^{h}-u\|^{2}_{1} \leq 2\|u^{h}-Q_{h}u\|^{2}_{1} + 2
\|Q_{h}u-u\|^{2}_{1} \leq ch^{2}\|u\|_{2}^{2}+ ch^{2} \leq ch^{2}.
\]
We conclude then with the desired error estimate.
\end{proof}

\section{Conclusion}

In this paper a dual mesh numerical scheme was proposed for a
nonlocal thermistor problem. We have showed the existence and
uniqueness of the approximate solution via Banach's f\/ixed point
theorem. We have also proved $H^1$-error bounds under minimal
regularity assumptions. We only obtain f\/irst-order estimates:
higher order estimates are dif\/f\/icult to obtain due to the
nonstandard nonlocal term. Optimal error analysis to the present
context, under appropriate smoothness assumptions on data, can be
derived by application of the techniques of \cite{clt}, but this
needs further developments.

\subsection*{Acknowledgements}

The support of the Portuguese Foundation for Science and
Technology (FCT) and post-doc fellowship SFRH/BPD/20934/2004 are
gratefully acknowledged. We would like to thank two anonymous
referees for valuable comments and suggestions.

\LastPageEnding

\end{document}